\documentclass[12pt]{article}
\usepackage{enumerate}
\usepackage{amssymb,a4wide,latexsym,makeidx,epsfig,fleqn}
\usepackage{amsthm}
\usepackage{amsmath}
\usepackage{enumerate}
\usepackage{float}
\newtheorem{theorem}{Theorem}[section]
\newtheorem{remark}[theorem]{Remark}
\newtheorem{example}[theorem]{Example}

\newtheorem{lemma}[theorem]{Lemma}

\newtheorem{corollary}[theorem]{Corollary}

\begin{document}
\textwidth 150mm \textheight 225mm
\title{Some upper bounds for the signless Laplacian spectral radius of digraphs
\thanks{ Supported by the National Natural Science
Foundation of China (No.11171273).}}
\author{{Weige Xi and Ligong Wang\footnote{Corresponding author.} }\\
{\small Department of Applied Mathematics, School of Science,}\\
{\small Northwestern Polytechnical University, Xi'an, Shaanxi 710072, P.R.China}
 \\{\small E-mail: xiyanxwg@163.com, lgwangmath@163.com }\\}
\date{}
\maketitle
\begin{center}
\begin{minipage}{120mm}
\vskip 0.3cm
\begin{center}
{\small {\bf Abstract}}
\end{center}
{\small Let $G=(V(G) ,E(G))$ be a digraph without loops and
multiarcs, where $V(G)=\{v_1,v_2,\ldots,v_n\}$ and $E(G)$ are the
vertex set and the arc set of $G$, respectively. Let $d_i^{+}$ be the
outdegree of the vertex $v_i$. Let $A(G)$ be the adjacency matrix of
$G$ and $D(G)=\textrm{diag}(d_1^{+},d_2^{+},\ldots,d_n^{+})$ be the
diagonal matrix with outdegrees of the vertices of $G$. Then we call
$Q(G)=D(G)+A(G)$ the signless Laplacian matrix of $G$. The spectral
radius of $Q(G)$ is called the signless Laplacian spectral radius of
$G$, denoted by $q(G)$. In this paper, some upper bounds for $q(G)$
are obtained. Furthermore, some upper bounds on
$q(G)$ involving outdegrees and the average 2-outdegrees of the
vertices of $G$ are also derived.

}
\end{minipage}
\end{center}

\section{Introduction }
\label{sec:ch6-introduction}

Let $G=(V(G), E(G))$ be a digraph without loops and multiarcs, where
$V(G)=\{v_1,v_2,$ $\ldots, v_n\}$ and $E(G)$ are the vertex set and the
arc set of $G$, respectively. If $(v_i, v_j)$ be an arc of $G$, then
$v_i$ is called the initial vertex of this arc and $v_j$ is called
the terminal vertex of this arc. For any vertex $v_i$ of $G$, we
denote $N_i^{+}=N_{v_i}^{+}(G)=\{v_j : (v_i,v_j)\in E(G) \}$ and
$N_i^{-}=N_{v_i}^{-}(G)=\{v_j : (v_j,v_i)\in E(G) \}$ the set of
out-neighbors and in-neighbors of $v_i$, respectively. Let
$d_i^{+}=|N_i^{+}|$ denote the outdegree of the vertex $v_i$ and
$d_i^{-}=|N_i^{-}|$ denote the indegree of the vertex $v_i$ in the
digraph $G$. The maximum vertex outdegree is denoted by $\Delta^+$,
and the minimum outdegree by $\delta^+$. If $\delta^+=\Delta^+$, then $G$ is a regular
digraph. Let $t_i^+=\sum\limits_{v_j\in N_i^{+}}d_j^{+}$ be the
2-outdegree of the vertex $v_i$, $m_i^{+}=\frac{t_i^{+}}{d_i^{+}}$
the average 2-outdegree of the vertex $v_i$. A digraph is strongly
connected if for every pair of vertices $v_i,v_j\in V(G)$, there
exists a directed path from $v_i$ to $v_j$ and a directed path from
$v_j$ to $v_i$. In this paper, we consider finite digraphs without
loops and multiarcs, which have at least one arc.

For a digraph $G$, let $A(G)=(a_{ij})$ denote the adjacency matrix
of $G$, where $a_{ij}=1$ if $(v_i,v_j)\in E(G)$ and $a_{ij}=0$ otherwise.
Let $D(G)=\textrm{diag}(d_1^{+},d_2^{+},\ldots,d_n^{+})$ be the
diagonal matrix with outdegrees of the vertices of $G$ and
$Q(G)=D(G)+A(G)$ the signless Laplacian matrix of $G$. However, the
signless Laplacian matrix of an undirected graph $D$ can be treated as
the signless Laplacian matrix of the digraph $G'$, where $G'$ is
obtained from $D$ by replace each edge with pair of oppositely
directed arcs joining the same pair of vertices. Therefore,
the research of the signless Laplacian matrix of a digraph has
more universal significance than undirected graph.

The eigenvalues of $Q(G)$ are called the signless Laplacian eigenvalues
of $G$, denoted by $q_1,q_2,\ldots,q_n$. In general $Q(G)$ are not
symmetric and so its eigenvalues can be complex numbers. We usually
assume that $|q_1|\geq|q_2|\geq\ldots\geq|q_n|$. The signless
Laplacian spectral radius of $G$ is denoted and defined as
$q(G)=|q_1|$, i.e., the largest absolute value of the signless
Laplacian eigenvalues of $G$. Since $Q(G)$ is a nonnegative matrix,
it follows from Perron Frobenius Theorem that $q(G)=q_1$ is a real
number.

For the  Laplacian spectral radius and signless Laplacian spectral radius of an undirected graph are
well treated in the literature, see
\cite{Wang,WYL,WeLi,Zhu} and \cite{ChWa,CTG,GDC,HaLu,HJZ,WeLi}, respectively. Recently, there are
some papers that give some lower or upper bounds for the spectral
radius of a digraph, see \cite{RB,GDa,XuXu}. Now we consider
the signless Laplacian spectral radius of a digraph $G$. For application it is
crucial to be able to computer or at least estimate $q(G)$ for a given digraph.

In 2013, S.B. Bozkurt and D. Bozkurt in \cite{BoBo} obtained the
following bounds for signless Laplacian spectral radius of a
digraph.
 \noindent\begin{equation}\label{eq:c1}
 q(G) \leq \max\{d_i^{+}+d_j^{+}: (v_i, v_j)\in E(G)\}.
\end{equation}
\noindent\begin{equation}\label{eq:c2}
 q(G) \leq \max\{d_i^{+}+m_i^{+}: v_i\in V(G)\}.
\end{equation}
\noindent\begin{equation}\label{eq:c3} q(G)\leq
\max\bigg\{\frac{d_i^{+}+d_j^{+}+
\sqrt{(d_i^{+}-d_j^{+})^{2}+4m_i^{+}m_j^{+}}}{2}: (v_i, v_j)\in
E(G)\bigg\}.
\end{equation}
\noindent\begin{equation}\label{eq:c4} q(G)\leq
\max\bigg\{d_i^{+}+\sqrt{\sum _{v_j:(v_j, v_i)\in E(G)}d_j^{+}}: v_i\in
V(G)\bigg\}.
\end{equation}

In 2014, Hong and You in \cite{HoYo} gave a sharp bound for the
signless Laplacian spectral radius of a digraph:
\noindent\begin{equation}\label{eq:c5} q(G)\leq \min_{1 \leq i \leq
n}\bigg\{\frac{d_1^{+}+2d_i^{+}-1+
\sqrt{(2d_i^{+}-d_1^{+}+1)^{2}+8\sum\limits_{k=1}^{i-1}(d_k^{+}-d_i^{+})}}{2}\bigg\}.
\end{equation}

\noindent\begin{remark}\label{re:c1} Note that $G$ is a strongly
connected digraph for bounds \eqref{eq:c1}, \eqref{eq:c3},
\eqref{eq:c4}, respectively.
\end{remark}

In this paper, we study on the signless Laplacian spectral radius of
a digraph $G$. We obtain some upper bounds for
$q(G)$, and we also show that some upper bounds on $q(G)$ involving
outdegrees and the average 2-outdegrees of the vertices of $G$ can
be obtained from our bounds.

\section{Preliminaries Lemmas}
\label{sec:1}

In this section, we give the following lemmas which will be used in
the following study.

\noindent\begin{lemma}\label{le:c1} (\cite{HoJo}) \ Let $M=(m_{ij})$
be an $n \times n$ nonnegative matrix with spectral radius $\rho(M)$, i.e.,
the largest eigenvalues of $M$,
and let $R_i=R_i(M)$ be the $i$-th row sum of $M$, i.e.,
$R_i(M)=\sum\limits_{j=1}^n m_{ij} \ (1 \leq i\leq n)$. Then
\begin{equation}\label{eq:ca}
\min\{R_{i}(M):1 \leq i\leq n\}\leq \rho(M)\leq \max\{R_{i}(M):1
\leq i\leq n\}.\end{equation}
 Moreover, if $M$ is irreducible, then any equality holds in \eqref{eq:ca} if
and only if $R_1=R_2=\ldots=R_n$.
\end{lemma}

\noindent\begin{lemma}\label{le:c2} (\cite{HoJo}) \ Let $M$ be an
irreducible nonnegative matrix. Then $\rho(M)$ is an eigenvalue of
$M$ and there is a positive vector $X$ such that $MX=\rho(M)X$.
\end{lemma}

\noindent\begin{lemma}\label{le:c6} (\cite{Li}) \ Let $A=(a_{ij})\in
\mathbb{C}^{n\times n}$, $r_i=\sum\limits_{j\neq i}|a_{ij}|$
for each $i=1,2,\ldots,n$,
$S_{ij}=\{z \in \mathbb{C}: |z-a_{ii}|\cdot|z-a_{jj}|\leq r_ir_j\}$
for all $i\neq j$. Also let $E(A)=\{(i,j) :a_{ij}\neq 0,1 \leq{i\neq j}\leq n\}$.
If $A$ is irreducible, then all eigenvalues of $A$ are contained in the
following region
\begin{equation}\label{eq:ce}
\Omega(A)=\bigcup_{(i,j) \in E(A)}S_{ij}.
\end{equation}
Furthermore, a boundary point $\lambda$ of \eqref{eq:ce} can be an
eigenvalue of $A$ only if  $\lambda$ locates on the boundary of each
oval region $S_{ij}$ for $e_{ij} \in E(A)$.
\end{lemma}

\section{ Some upper bounds for the signless Laplacian
spectral radius of digraphs}
\label{sec:2}

 In this section, we present some upper bounds for the signless
Laplacian spectral radius $q(G)$ of a digraph $G$
and also show that some bounds involving outdegrees, the average
2-outdegrees, the maximum outdegree and the minimum outdegree of the
vertices of $G$ with $n$ vertices and $m$ arcs can be obtained from
our bounds.

\noindent\begin{theorem}\label{th:c7}  \ Let $G$ be a strongly
connected digraph with $n\geq 3$ vertices, $m$ arcs, the maximum
vertex outdegree $\Delta^+$ and the minimum outdegree $\delta^{+}$.
Then
\begin{equation}\label{eq:c27}
q(G)\leq
\max\{\Delta^{+}+\delta^{+}-1+\frac{m-\delta^{+}(n-1)}{\Delta^{+}},
\delta^{+}+1+\frac{m-\delta^{+}(n-1)}{2} \}.
\end{equation} Moreover, if
$G(\neq\overset{\longrightarrow}{C_{n}})$ is a regular digraph or
$G\cong \overset{\longleftrightarrow}{K}_{1,n-1}$, where
$\overset{\longleftrightarrow}{K}_{1,n-1}$ denotes the digraph on
$n$ vertices which replace each edge in star graph $K_{1,n-1}$ with
the pair of oppositely directed arcs, then the equality holds in
\eqref{eq:c27}
\end{theorem}

\begin{proof}  \ From \eqref{eq:c2}, we know that
$q(G)\leq \max\{d_i^{+}+m_i^{+}: v_i\in V(G)\}.$ So we only need to
prove that $\max\{d_i^{+}+m_i^{+}: v_i\in V(G)\}\leq
\max\{\Delta^{+}+\delta^{+}-1+\frac{m-\delta^{+}(n-1)}
{\Delta^{+}},\delta^{+}+1+\frac{m-\delta^{+}(n-1)}{2} \}$. Suppose
$\max\{d_i^{+}+m_i^{+}: v_i\in V(G)\}$ occurs at vertex $u$. Two
cases arise $d_u^{+}=1$, or $2\leq d_u^{+}\leq \Delta^{+}$.

\noindent{\textbf{{Case 1.}}} \ $d_u^{+}=1.$ Suppose that
$N_u^{+}=\{w\}.$ Since $m_u^{+}=d_w^{+}\leq \Delta^{+},$ thus
$d_u^{+}+m_u^{+}\leq 1+\Delta^{+}$. Since $\sum\limits_{v_i\in
V(G)}d_i^{+}=m,$ let $d_j^{+}=\Delta^{+}$, then $\sum\limits_{i\neq
j}d_i^{+}=m- \Delta^{+}\geq (n-1)\delta^{+}$, so
$m-(n-1)\delta^{+}\geq \Delta^{+}$.
 Therefore $\delta^{+}-1+\frac{m-\delta^{+}(n-1)}{\Delta^{+}}\geq
\delta^{+}-1+\frac{\Delta^{+}}{\Delta^{+}}=\delta^{+}\geq 1$. Thus
$d_u^{+}+m_u^{+}\leq 1+\Delta^{+}\leq
\Delta^{+}+\delta^{+}-1+\frac{m-\delta^{+}(n-1)}{\Delta^{+}}$, the
result follows.

\noindent{\textbf{{Case 2.}}} \ $2\leq d_u^{+}\leq \Delta^{+}$. Note
that $m-(n-1)\delta^{+}\geq d_u^{+}\geq 2$, and
\begin{eqnarray*}
m&=&\sum\limits_{v:(u,v)\in E(G)}d_v^{+}+\sum\limits_{
v:(u,v)\notin E(G)}d_v^{+}\\
&\geq &\sum\limits_{v:(u,v)\in E(G)}d_v^{+}+
d_u^{+}+(n-d_u^{+}-1)\delta^{+},
\end{eqnarray*}
 thus
\begin{eqnarray*}
\sum\limits_{v:(u,v)\in E(G)}d_v^{+}&\leq&
m-d_u^{+}-(n-d_u^{+}-1)\delta^{+}\\
&=&m-(n-1)\delta^{+}+(\delta^{+}-1)d_u^{+}
\end{eqnarray*}
$$m_u^{+}=\frac{\sum\limits_{v:(u,v)\in E(G)}d_v^{+}}{d_u^{+}}\leq
\frac{m-(n-1)\delta^{+}}{d_u^{+}}+\delta^{+}-1.$$ This follows that
$m_u^{+}+d_u^{+}\leq
d_u^{+}+\frac{m-(n-1)\delta^{+}}{d_u^{+}}+\delta^{+}-1$. Let
$f(x)=x+\frac{m-(n-1)\delta^{+}}{x}+\delta^{+}-1$, where $x\in
[2,\Delta^{+}].$ It is easy to see that
$f'(x)=1-\frac{m-(n-1)\delta^{+}}{x^{2}}$. Let
$a=m-(n-1)\delta^{+},$ then $\sqrt{a}$ is the unique positive root
of $f'(x)=0.$ We consider the next three Subcases.

\noindent{\textbf{{Subcase 1.}}} \ $\sqrt{a}< 2.$ When $x\in
[2,\Delta^{+}],$ since $f'(x)>0,$ then  $f(x)\leq f(\Delta^{+})$.

\noindent{\textbf{{Subcase 2.}}} \ $2\leq\sqrt{a}\leq \Delta^{+}.$
Then $f'(x)<0$ for $x\in [2,\sqrt{a})$, and $f'(x)\geq 0$, for $x\in
[\sqrt{a}, \Delta^{+}]$. Thus, $f(x)\leq \max\{f(2),
f(\Delta^{+})\}$.

\noindent{\textbf{{Subcase 3.}}} \ $\Delta^{+}<\sqrt{a}$. When $x\in
[2,\Delta^{+}]$, since $f'(x)<0$, then $f(x)\leq f(2)$.

Recall that $2\leq d_u^{+}\leq \Delta^{+}$, thus
$$ m_u^{+}+d_u^{+}\leq \max\{f(2), f(\Delta^{+})\}$$
$$=\max\{\Delta^{+}+\delta^{+}-1+\frac{m-\delta^{+}(n-1)}{\Delta^{+}},
\delta^{+}+1+\frac{m-\delta^{+}(n-1)}{2}\}.$$ If
$G(\neq\overset{\longrightarrow}{C_{n}})$  is a regular digraph,
then $d_i^{+}+m_i^{+}=2d_i^{+}=2\Delta^{+}$ for all $v_i\in V(G)$.
We can get $q(G)=2\Delta^{+}$. Since
$G(\neq\overset{\longrightarrow}{C_{n}})$ is a strongly connected digraph,
then we may assume that $\Delta^{+} \geq 2$, this implies that
$\delta^{+}+1+\frac{m-\delta^{+}(n-1)}{2}=\Delta^{+}+1+\frac{\Delta^{+}}{2}\leq
2\Delta^{+}=\Delta^{+}+\delta^{+}-1+\frac{m-\delta^{+}(n-1)}{\Delta^{+}}$.
So
$\max\{\Delta^{+}+\delta^{+}-1+\frac{m-\delta^{+}(n-1)}{\Delta^{+}},
\delta^{+}+1+\frac{m-\delta^{+}(n-1)}{2} \}=2\Delta^{+}$. Thus, the
equality holds. If $G\cong \overset{\longleftrightarrow}{K}_{1,n-1}
$, we can get $q(G)=n$. Since
$\delta^{+}+1+\frac{m-\delta^{+}(n-1)}{2}=2+\frac{n-1}{2}\leq n$
from $n\geq 3$ and
$\Delta^{+}+\delta^{+}-1+\frac{m-\delta^{+}(n-1)}{\Delta^{+}}=n-1+1-1+\frac{n-1}{n-1}=n$.
So
$\max\{\Delta^{+}+\delta^{+}-1+\frac{m-\delta^{+}(n-1)}{\Delta^{+}},
\delta^{+}+1+\frac{m-\delta^{+}(n-1)}{2} \}=n$. Thus, the equality
also holds. By combining the above discussion, the result follows.
\end{proof}

\noindent\begin{corollary}\label{co:ck}
 Let $G$ be a strongly
connected digraph with $n\geq 3$ vertices, $m$ arcs, the maximum
 outdegree $\Delta^+$ and the minimum outdegree $\delta^{+}$.
If $\Delta^{+}\geq\frac{m-(n-1)}{2}$ and $\delta^{+}=1$, then
\begin{equation}\label{eq:cm}
q(G)\leq \Delta^{+}+2.
\end{equation}
\end{corollary}

\begin{proof} Because $\Delta^{+}+\delta^{+}-1+\frac{m-\delta^{+}(n-1)}{\Delta^{+}}\leq \Delta^{+}+2$,
$\delta^{+}+1+\frac{m-\delta^{+}(n-1)}{2}\leq \Delta^{+}+2$, therefore by Theorem \ref{th:c7}, we have
$q(G)\leq \Delta^{+}+2.$
\end{proof}

Let $G^{*}(m,n,\frac{m-(n-1)}{2},1)$ be a class of strongly
connected digraphs with $\Delta^{+}\geq\frac{m-(n-1)}{2}$, $\delta^{+}=1$,
and there exists a vertex $v_0\in V(G)$ such that $d_{v_0}=\Delta^{+}$ and there exists
a vertex $v_k\in N_{v_0}^{+}$, $d_{v_k}^{+}\geq 2$.

\noindent\begin{remark}\label{re:c2}
For $G\in G^{*}(m,n,\frac{m-(n-1)}{2},1)$, we have
$\Delta^{+}+2\leq\max\{d_i^{+}+d_j^{+}: (v_i, v_j)\in E(G)\}$, thus the upper bound
\eqref{eq:cm} is better than the upper bound \eqref{eq:c1} for the class of digraphs
$G\in G^{*}(m,n,\frac{m-(n-1)}{2},1)$. But for general digraphs, the upper bound
\eqref{eq:cm} is incomparable with the upper bound \eqref{eq:c1}.
\end{remark}
\begin{example}
Let $G$ be the digraph of order 4, as shown in Figure 1. Since it has 9 arcs, and
the maximum outdegree $\Delta^+=3=\frac{9-(4-1)}{2}$, the minimum outdegree $\delta^{+}=1$,
and there exists a vertex $v_4\in N_{v_1}^{+}$, $d_{v_4}^{+}=3>2$,
therefore $G=G^{*}(9,4,3,1)$.
\begin{figure}[H]
\begin{centering}
%\subfigure[$G_1$] {\label{G_0}
\includegraphics[scale=0.5]{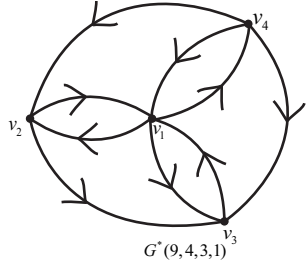}
 \caption{Graph $G^*(9,4,3,1)$}
\end{centering}
\end{figure}
\end{example}
\begin{table}[H]
\centering\caption{Values of the upper bounds for example 1.}

\begin{tabular}{cccc}
\hline
&$q(G)$&\eqref{eq:c1}&\eqref{eq:cm}\\
 \hline
$G^{*}(9,4,3,1)$ &  4.7321& 6&5  \\
\hline
\end{tabular}
\end{table}

\noindent\begin{theorem}\label{th:c8}  \ Let $G$ be a strongly
connected digraph  with vertex set $V(G)=\{v_1,v_2,$ $\ldots, v_n\}$
and arc set $E(G)$. Then
\begin{equation}\label{eq:c28}
 q(G)\leq \max\bigg\{\frac{d_i^{+}+d_j^{+}+
\sqrt{(d_i^{+}-d_j^{+})^{2}+{4\sqrt{d_i^{+}m_i^{+}}\sqrt{d_j^{+}m_j^{+}}}}}{2}
:(v_i, v_j)\in E(G) \bigg\}.\end{equation} Moreover if $G$ is a
regular digraph or a bipartite semiregular digraph, then the equality holds in \eqref{eq:c28}.
\end{theorem}

\begin{proof}  \ From the definition of $D=D(G)$ we get
$D^{\frac{1}{2}}=\textrm{diag}(\sqrt{d_i^{+}}:v_i\in V(G)),$ and
consider the similar matrix $P=D^{-\frac{1}{2}}Q(G)D^{\frac{1}{2}}.$
Since $G$ is a strongly connected digraph, it is easy to see that
$P$ is irreducible and nonnegative. Now the $(i,j)$-th element of
$P=D^{-\frac{1}{2}}Q(G)D^{\frac{1}{2}}$ is
$$p_{ij} =
\begin{cases}
\ d_i^{+}       & \textrm{ if $i=j $,}\\
\ \frac{\sqrt{d_j^{+}}}{\sqrt{d_i^{+}}}  & \textrm{ if $(v_i, v_j)\in E(G) $,}\\
\ 0               &  \textrm{ otherwise.}
\end{cases}$$
Let $R_i(P)$ be the $i$-th row sum of $P$ and
$R_i^{'}(P)=R_i-d_i^{+}.$ Then by Cauchy-Schwarz inequality, we have
\begin{align*}
R_i^{'}(P)^{2}=&\left(\sum\limits_{v_j:(v_i,v_j)\in
E(G)}\frac{\sqrt{d_j^{+}}}{\sqrt{d_i^{+}}} \right)^{2} \leq
\sum\limits_{v_j:(v_i,v_j)\in E(G)}1^{2}\sum\limits_{v_j:(v_i,v_j)\in
E(G)}\frac{d_j^{+}}{d_i^{+}}\\
 =&\sum\limits_{v_j:(v_i,v_j)\in
E(G)}d_j^{+} =d_i^{+}m_i^{+}.
\end{align*}
 Since $P$ is irreducible and
nonnegative, $\rho(P)$ denotes the spectral radius of $P$. Then by
Lemma \ref{le:c6}, there at least exists $(v_i,v_j)\in E(G)$ such
that $\rho(P)$ is contained in the following oval region
$$|\rho(P)-d_i^{+}||\rho(P)-d_j^{+}|\leq R_i^{'}(P)R_i^{'}(P)
\leq\sqrt{d_i^{+}m_i^{+}}\sqrt{d_j^{+}m_j^{+}}.$$ Obviously,
$\rho(P)=q(G)>\max\{d_i^{+}: v_i \in E(G)\},$ and
$(\rho(P)-d_i^{+})(\rho(P)-d_j^{+})\leq|\rho(P)-d_i^{+}||\rho(P)-d_j^{+}|.$
Therefore, solving the above inequality we obtain
$$ q(G)\leq\bigg\{\frac{d_i^{+}+d_j^{+}+\sqrt{(d_i^{+}-d_j^{+})^{2}
+{4\sqrt{d_i^{+}m_i^{+}}\sqrt{d_j^{+}m_j^{+}}}}}{2} \bigg\}.$$ Hence
\eqref{eq:c28} holds.

If $G$ is a regular digraph,
$$q(G)=2\Delta^{+}=\max\bigg\{\frac{d_i^{+}+d_j^{+}+\sqrt{{(d_i^{+}-d_j^{+})^{2}
+4\sqrt{d_i^{+}m_i^{+}}\sqrt{d_j^{+}m_j^{+}}}}}{2} :(v_i, v_j)\in
E(G) \bigg\}.$$ Thus, the equality holds.

If $G$ is a bipartite semiregular digraph,
$$\max\bigg\{\frac{d_i^{+}+d_j^{+}+\sqrt{{(d_i^{+}-d_j^{+})^{2}
+4\sqrt{d_i^{+}m_i^{+}}\sqrt{d_j^{+}m_j^{+}}}}}{2} :(v_i, v_j)\in
E(G) \bigg\}=d_i^{+}+d_j^{+}.$$ Because $q(G)=\rho(D^{-1}Q(G)D)$, the $i$-th row sum of
$D^{-1}Q(G)D$ is $d_i^{+}+m_i^{+}$, and $G$ is a bipartite semiregular digraph, therefore
$d_i^{+}+m_i^{+}=d_i^{+}+d_j^{+}, (v_i, v_j)\in E(G)$, that is the row sums of $D^{-1}Q(G)D$
are all equal, then by Lemma \ref{le:c1}, $\rho(D^{-1}Q(G)D)=d_i^{+}+d_j^{+}$.
Thus we have
\begin{align*}
q(G)&=\rho(D^{-1}Q(G)D)=d_i^{+}+d_j^{+} \\
&=\max\bigg\{\frac{d_i^{+}+d_j^{+}+\sqrt{{(d_i^{+}-d_j^{+})^{2}
+4\sqrt{d_i^{+}m_i^{+}}\sqrt{d_j^{+}m_j^{+}}}}}{2} :(v_i, v_j)\in
E(G) \bigg\}.
\end{align*}
Then the equality holds.
\end{proof}

 For a digraph $G=(V(G), E(G)),$ let $f : V(G)\times V(G)
\rightarrow \mathbb{R}$ be a function. If $f(v_i, v_j)>0$ for all
$(v_i, v_j)\in E(G),$ we say $f$ is positive on arcs.

\noindent\begin{theorem}\label{th:c9} \ Let $G=(V(G), E(G))$ be a
digraph. Let $f : V(G)\times V(G)\rightarrow \mathbb{R^{+}}\bigcup
\{0\}$ be a nonnegative function which is positive on arcs. Then
\begin{equation}\label{eq:c29}
q(G)\leq \max\left\{\frac{\sum\limits_{v_k:(v_i, v_k)\in E(G)}f(v_i,
v_k)+ \sum\limits_{v_k:(v_j, v_k)\in E(G)}f(v_j, v_k)}{f(v_i,v_j)} :
(v_i,v_j)\in E(G)\right\}. \end{equation}
\end{theorem}

\begin{proof} \ Let ${\bf X}=(x_1,x_2,\ldots,x_n)^{T}$ be
an eigenvector corresponding to the eigenvalue $q(G)$ of $Q(G)$.
Since
$$Q(G){\bf X}=q(G){\bf X}.$$
Then for $1 \leq i\leq n$
\begin{equation}\label{eq:c30}
q(G)x_i=d_i^{+}x_i+\sum\limits_{v_k:(v_i,v_k)\in
E(G)}x_k=\sum\limits_{v_k:(v_i,v_k)\in E(G)}(x_i+x_k).
\end{equation}
By \eqref{eq:c30}, we have
$$q(G)(x_i+x_j)=\sum\limits_{v_k:(v_i,v_k)\in E(G)}(x_i+x_k)+\sum\limits_{v_k:(v_j,v_k)\in E(G)}(x_j+x_k).$$
For convenience we use $f(i, j)$ denote $f(v_i, v_j).$ Set $g(i,
j)=\frac{x_i+x_j}{f(i,j)}.$ If $(v_i, v_j)\in E(G),$ then
\begin{equation}\label{eq:c31}
q(G)f(i,j)g(i,j)=\sum\limits_{v_k:(v_i,v_k)\in
E(G)}f(i,k)g(i,k)+\sum\limits_{v_k:(v_j,v_k)\in E(G))}f(j,k)g(j,k).
\end{equation}
By \eqref{eq:c31}, we get
\begin{align*}
|q(G)f(i,j)g(i,j)|=&q(G)f(i,j)|g(i,j)| \\
\leq&\sum\limits_{v_k:(v_i,v_k)\in
E(G)}f(i,k)|g(i,k)|+\sum\limits_{v_k:(v_j,v_k)\in E(G)}f(j,k)|g(j,k)|.
\end{align*} Now choose $i_1,j_1$ such that $(v_{i_1}, v_{j_1})\in
E(G)$ and $|g(i_1, j_1)|=\max\{|g(i, j)| : (v_i, v_j)\in E(G)\}.$ If
$g|(i_1, j_1)|=0,$ then $|g(i, j)|=0 $ for all arcs $(v_i, v_j)\in
E(G).$ i.e.,$x_i+x_j=0$ for all arcs $(v_i, v_j)\in E(G).$ By
\eqref{eq:c30}, we have $q(G)=0$ which is impossible, since $G$ has
at least one arc. So $|g(i_1,j_1)|>0.$ Then
$$q(G)f(i_1,j_1)|g(i_1,j_1)|\leq\sum\limits_{v_k:(v_{i_1},v_k)\in E(G)}f(i_1,k)|g(i_1,k)|+
\sum\limits_{v_k:(v_{j_1},v_k)\in E(G)}f(j_1,k)|g(j_1,k)|.$$ Therefore,
we obtain
\begin{align*}
q(G)\leq&\sum\limits_{v_k:(v_{i_1},v_k)\in
E(G)}\frac{f(i_1,k)}{f(i_1,j_1)}\frac{|g(i_1,k)|}{|g(i_1,j_1)|}+
\sum\limits_{v_k:(v_{j_1},v_k)\in E(G)}\frac{f(j_1,k)}{f(i_1,j_1)}\frac{|g(j_1,k)|}{|g(i_1,j_1)|}\\
    \leq&\sum\limits_{v_k:(v_{i_1},v_k)\in E(G)}\frac{f(i_1,k)}{f(i_1,j_1)}+
    \sum\limits_{v_k:(v_{j_1},v_k)\in E(G)}\frac{f(j_1,k)}{f(i_1,j_1)},
\end{align*}
i.e.,$$q(G)\leq\frac{\sum\limits_{v_k:(v_{i_1},v_k)\in
E(G)}f(i_1,k)+\sum\limits_{v_k:(v_{j_1},v_k)\in
E(G)}f(j_1,k)}{f(i_1,j_1)}, \textrm {where $(v_{i_1}, v_{j_1})\in
E(G)$}.$$ This proves the desired result.
\end{proof}

\noindent\begin{corollary}\label{co:c10} \ Let $G=(V(G), E(G))$ be a
digraph. Then
\begin{equation}\label{eq:c32}q(G)\leq \max\left\{d_i^{+}
\sqrt{\frac{m_i^{+}}{d_j^{+}}}+d_j^{+}\sqrt{\frac{m_j^{+}}{d_i^{+}}}
: {(v_i, v_j)\in E(G)}\right\}. \end{equation}
\end{corollary}

\begin{proof} \ Setting $f(v_i,
v_j)=\sqrt{d_i^{+}d_j^{+}}$ in \eqref{eq:c29}, by Cauchy-Schwarz
inequality, $$\sum\limits_{v_k:(v_i, v_k)\in E(G)}f(v_i, v_k)
=\sum\limits_{v_k:(v_i, v_k)\in
E(G)}\sqrt{d_i^{+}d_k^{+}}=\sum\limits_{v_k:(v_i, v_k)\in
E(G)}(\sqrt{d_i^{+}}\sqrt{d_k^{+}})$$ $$\leq\sqrt{\sum\limits_{v_k:(v_i,
v_k)\in E(G)}d_i^{+}\sum\limits_{v_k:(v_i, v_k)\in E(G)}d_k^{+}}=
\sqrt{{d_i^{+}}^{2}\sum\limits_{v_k:(v_i, v_k)\in
E(G)}d_k^{+}}=d_i^{+}\sqrt{d_i^{+}m_i^{+}}.$$ By \eqref{eq:c29}, we
get
\begin{align*}
 q(G)\leq& \max\left\{\frac{\sum\limits_{v_k:(v_i, v_k)\in E(G)}f(v_i, v_k)+
\sum\limits_{v_k:(v_j, v_k)\in E(G)}f(v_j, v_k)}{f(v_i,v_j)} :
(v_i,v_j)\in E(G)\right\} \\
     \leq& \max\left\{\frac{d_i^{+}\sqrt{d_i^{+}m_i^{+}}+d_j^{+}\sqrt{d_j^{+}m_j^{+}}}
     {\sqrt{d_i^{+}d_j^{+}}} : (v_i,v_j)\in E(G)\right\} \\
     =&\max\left\{d_i^{+}\sqrt{\frac{m_i^{+}}{d_j^{+}}}+d_j^{+}\sqrt{\frac{m_j^{+}}{d_i^{+}}}
: {(v_i, v_j)\in E(G)}\right\}.
\end{align*}
\end{proof}

\noindent\begin{corollary}\label{co:c13} \  Let $G=(V(G), E(G))$ be
a digraph. Then
\begin{equation}\label{eq:c35}q(G)\leq
\max\left\{\frac{d_i^{+}(d_i^{+}+m_i^{+})+d_j^{+}(d_j^{+}+m_j^{+})}{d_i^{+}+d_j^{+}}
: (v_i,v_j)\in E(G)\right\}. \end{equation}
\end{corollary}

\begin{proof} \ Setting $f(v_i, v_j)=d_i^{+}+d_j^{+}$ in
\eqref{eq:c29}, since $\sum\limits_{v_k:(v_i, v_k)\in E(G)}f(v_i,
v_k)=\sum\limits_{v_k:(v_i, v_k)\in
E(G)}(d_i^{+}+d_k^{+})=d_i^{+}(d_i^{+}+m_i^{+})$, So we get the
desired result.
 \end{proof}

\noindent\begin{corollary}\label{co:c11} \ Let $G=(V(G), E(G))$ be a
digraph. Then
\begin{equation}\label{eq:c33}q(G)\leq \max\left\{\frac{d_i^{+}\sqrt{d_i^{+}
+m_i^{+}}+d_j^{+}\sqrt{d_j^{+}+m_j^{+}}}{\sqrt{d_i^{+}+d_j^{+}}}
:{(v_i, v_j)\in E(G)}\right\}.
\end{equation}
\end{corollary}

\begin{proof} \ Setting $f(v_i,
v_j)=\sqrt{d_i^{+}+d_j^{+}}$ in \eqref{eq:c29}, since

$\sum\limits_{v_k:(v_i, v_k)\in E(G)}f(v_i, v_k)=\sum\limits_{v_k:(v_i,
v_k)\in E(G)}(1\cdot
\sqrt{d_i^{+}+d_k^{+}})\leq\sqrt{d_i^{+}\sum\limits_{v_k:(v_i, v_k)\in
E(G)}(d_i^{+}+d_k^{+})}$

$=\sqrt{d_i^{+}({d_i^{+}}^{2}+d_i^{+}m_i^{+})}=d_i^{+}\sqrt{d_i^{+}+m_i^{+}}$
by Cauchy-Schwarz inequality.

Thus by \eqref{eq:c29} we get the
desired result.
 \end{proof}

\noindent\begin{corollary}\label{co:c12} \ Let $G=(V(G), E(G))$ be a
digraph. Then
\begin{equation}\label{eq:c34}q(G)\leq \max\left\{\frac{d_i^{+}
(\sqrt{d_i^{+}}+\sqrt{m_i^{+}})+d_j^{+}(\sqrt{d_j^{+}}+\sqrt{m_j^{+}})}
{\sqrt{d_i^{+}}+\sqrt{d_j^{+}}}:{(v_i, v_j)\in E(G)}\right\}.
\end{equation}
\end{corollary}

\begin{proof} \ Setting $f(v_i,
v_j)=\sqrt{d_i^{+}}+\sqrt{d_j^{+}}$ in \eqref{eq:c29}, since
$\sum\limits_{v_k:(v_i, v_k)\in E(G)}f(v_i, v_k)=\sum\limits_{v_k:(v_i,
v_k)\in E(G)}(\sqrt{d_i^{+}}$
$+\sqrt{d_k^{+}})=d_i^{+}\sqrt{d_i^{+}}+\sum\limits_{v_k:(v_i, v_k)\in
E(G)}(1\cdot \sqrt{d_k^{+}})\leq
{d_i^{+}}^{\frac{3}{2}}+\sqrt{d_i^{+}\sum\limits_{v_k:(v_i, v_k)\in
E(G)}d_k^{+}}=d_i^{+}(\sqrt{d_i^{+}}+\sqrt{m_i^{+}})$ by
Cauchy-Schwarz inequality. By \eqref{eq:c29} the result follows.
\end{proof}

Notice that \eqref{eq:c33} and \eqref{eq:c34} can be viewed as adding square roots to
\eqref{eq:c35} at difference places.

\section{Example}
\label{sec:3}

Let $G_1$, $G_2$ be the digraphs of order 4,6, respectively, as shown in Figure 2.
\begin{figure}[H]
\begin{centering}
%\subfigure[$G_1$] {\label{G_0}
\includegraphics[scale=1.5]{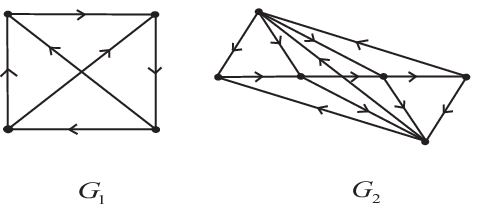}
 \caption{}
\end{centering}
\end{figure}

\begin{table}[H]
\centering\caption{Values of the various bounds for example 1.}
\begin{tabular}{cccccccccccccccc}
\hline
&$q(G)$&\eqref{eq:c1}&\eqref{eq:c2}&\eqref{eq:c3}&\eqref{eq:c4}&\eqref{eq:c5} \\
 &\eqref{eq:c27}&\eqref{eq:c28}&\eqref{eq:c32}&\eqref{eq:c35}&\eqref{eq:c33}&\eqref{eq:c34}\\
 \hline
$G_1$ &  3.0000& 4.0000&  3.5000&  3.3028 &  3.4142& 3.5616 \\
& 3.5000 & 3.5651 & 3.4495& 3.3333& 3.6029& 3.5731  \\
\hline
     $G_2$ &   4.1984& 5.0000&   4.6667 &  4.6016&  5.0000 &  4.7321  \\
 & 5.5000& 4.7913& 4.5644& 4.6000& 4.7956 & 4.7866 \\
\hline
\end{tabular}
\end{table}

\noindent\begin{remark}\label{re:c3}
Obviously, from Table 1, the bound \eqref{eq:c3} is the best in all known upper bounds for $G_1$, and the bound \eqref{eq:c32} is the best for $G_2$. Finally bound \eqref{eq:c35} is the second-best bounds for $G_1$ and $G_2$. In general, these bounds are incomparable.
\end{remark}
%\section*{Acknowledgements}

%The authors are very grateful to the referees for many detailed
%comments and suggestions, which are very helpful for improving the
%presentation of the manuscript. Lemma 2.8 and Theorem 3.9 are due
%to an anonymous referee, to whom we would like to express our
%special thanks.

\end{document}